\documentclass[12pt,leqno]{article}
\usepackage{amsmath, amsthm, amsfonts, amssymb, color}
\usepackage{mathrsfs}
\theoremstyle{definition}

\newcommand{\scr}[1]{\mathscr #1}
\definecolor{wco}{rgb}{0.5,0.2,0.3}

\numberwithin{equation}{section} \theoremstyle{remark}

\newcommand{\ua}{\uparrow}

\title{{\bf Harnack Inequalities on Manifolds with Boundary and Applications}\footnote{Supported in
 part by WIMCS, NNSFC(10721091) and the 973-Project.}
}
\author{
{\bf Feng-Yu Wang}\\
\footnotesize{School of Mathematical Sci. and Lab. Math. Com. Sys.,
Beijing Normal
University, Beijing 100875, China}\\
\footnotesize{and}\\ \footnotesize{Department of Mathematics,
Swansea University, Singleton Park, SA2 8PP, UK}\\ \footnotesize{Email: wangfy@bnu.edu.cn;
F.Y.Wang@swansea.ac.uk}}
\begin{document}
\def\R{\mathbb R}  \def\ff{\frac} \def\ss{\sqrt} \def\BB{\mathbb
B}
\def\N{\mathbb N} \def\kk{\kappa} \def\m{{\bf m}}
\def\dd{\delta} \def\DD{\Delta} \def\vv{\varepsilon} \def\rr{\rho}
\def\<{\langle} \def\>{\rangle} \def\GG{\Gamma} \def\gg{\gamma}
  \def\nn{\nabla} \def\pp{\partial} \def\tt{\tilde}
\def\d{\text{\rm{d}}} \def\bb{\beta} \def\aa{\alpha} \def\D{\scr D}
\def\EE{\mathbb E} \def\si{\sigma} \def\ess{\text{\rm{ess}}}
\def\beg{\begin} \def\beq{\begin{equation}}  \def\F{\scr F}
\def\Ric{\text{\rm{Ric}}} \def\Hess{\text{\rm{Hess}}}\def\B{\scr B}
\def\e{\text{\rm{e}}} \def\ua{\underline a} \def\OO{\Omega}  \def\oo{\omega}
 \def\tt{\tilde} \def\Ric{\text{\rm{Ric}}}
\def\cut{\text{\rm{cut}}} \def\P{\mathbb P} \def\ifn{I_n(f^{\bigotimes n})}
\def\C{\scr C}      \def\aaa{\mathbf{r}}     \def\r{r}
\def\gap{\text{\rm{gap}}} \def\prr{\pi_{{\bf m},\varrho}}  \def\r{\mathbf r}
\def\Z{\mathbb Z} \def\vrr{\varrho} \def\ll{\lambda}
\def\L{\scr L}\def\Tt{\tt} \def\TT{\tt}\def\II{\mathbb I}
\def\i{{\rm i}}\def\Sect{{\rm Sect}}\def\E{\mathbb E} \def\H{\mathbb H}
\def\M{\scr M}

\maketitle
\begin{abstract} On a large class of Riemannian manifolds with boundary, some dimension-free  Harnack inequalities
for the Neumann semigroup  is proved to be equivalent to the
convexity of the boundary and a curvature condition. In particular,
for $p_t(x,y)$   the Neumann heat kernel w.r.t. a volume type
measure $\mu$ and for $K$ a constant, the curvature condition
$\Ric-\nn Z\ge K$ together with the convexity of the boundary  is
equivalent to the heat kernel entropy inequality

$$\int_M p_t(x,z)\log \ff{p_t(x,z)}{p_t(y,z)}\,\mu(\d z)\le \ff{K\rr(x,y)^2}{2(\e^{2Kt}-1)} ,
\ \ t>0, x,y\in M,$$   where $\rr$ is the Riemannian distance. The
main result is partly extended to
  manifolds with non-convex boundary and applied to derive  the HWI inequality.

\end{abstract} \noindent
 \vskip 2cm

\section{Introduction} Let $M$ be a connected complete Riemannian manifold possibly with a boundary  $\pp M$.
Let $L=\DD +Z$ for a $C^2$ vector field $Z$ on $M$. Let $P_t$ be the
(Neumann if $\pp M\ne\emptyset$) diffusion semigroup generated by
$L$. Then for any measure $\mu$ equivalent to the Riemannian volume,
$P_t$ has a heat kernel $\{p_t(x,y): x,y\in M\}$ with respect to
$\mu$, i.e.

$$P_t f(x)= \int_M p_t(x,y) f(y)\mu(\d y)$$ holds for any bounded measurable function
$f$.  When $\pp M=\emptyset$, there exist many equivalent statements
on the semigroup $P_t$ for the following curvature condition (known
as the $\GG_2$ condition of Bakry and Emery \cite{BE}):

\beq\label{1.1} \Ric (X,X)-\<\nn_X Z,X\> \ge -K|X|^2,\ \ \ X\in
TM,\end{equation} where $K\in \R$ is a constant. See e.g. \cite{B,
BL} for equivalent gradient and Poincar\'e/log-Sobolev inequalities,
\cite{S} for equivalent cost (or Wasserstein distance) inequalities,
and \cite{W04} for equivalent dimension-free Harnack inequalities.
These equivalences also hold if $M$ has a convex boundary (cf.
\cite{W04}).  The main purpose of this paper is to provide
equivalent heat kernel inequalities for (\ref{1.1}) and the
convexity of $\pp M$. To this end we first recall two  known Harnack type
inequalities for $P_t$.

According to \cite[Lemma 2.2]{W97}, if $\pp M$ is either empty or
convex, then (\ref{1.1}) implies the Harnack inequality

\beq\label{Harnack} \ff{(P_t f(x))^\aa}{P_tf^\aa f(y)}\le
\exp\Big[\ff{K\aa\rr(x,y)^2}{2(\aa-1)(1-\e^{-2Kt})}\Big],\ \ f\in
\M_b^+(M),   t>0, x,y\in M\end{equation} for all $\aa>1$, where
$\M_b^+(M)$ is the set of all positive measurable functions on $M$, and
$\rr$ is the Riemannian distance on $M$. It is also proved in
\cite{W04} that, if (\ref{Harnack}) holds for all $\aa>1$ then
(\ref{1.1}) holds. In this paper we shall prove that (\ref{Harnack}) is equivalent to (\ref{1.1}) for each fixed $\aa>1.$

Next, when $\pp M$ is either empty or convex, we prove that
(\ref{1.1}) is also equivalent to the following log-Harnack
inequality, a limit version of (\ref{Harnack}) as $\aa\to\infty$
(see Section 2):

\beq\label{Log} P_t(\log f)(x)\le \log P_t f(y) +
\ff{K\rr(x,y)^2}{2(1-\e^{-2Kt})},\ \ \ f\ge 1, t>0, x,y\in
M.\end{equation} Note that this type inequality was used in
\cite{BGL} for the study of HWI inequalities on manifolds without boundary.
In conclusion we have the following result.

\beg{thm} \label{T1.1} Assume that $\pp M$ is either empty or
convex.  Let $K\in \R.$   Then the following statements are equivalent to each other:
\beg{enumerate} \item[$(1)$] $\Ric(X,X)-\<\nn_X Z,X\>\ge -K|X|^2,\ \ X\in TM.$
\item[$(2)$] The Harnack inequality
$(\ref{Harnack})$ holds for all $\aa>1.$ \item[$(3)$] The Harnack inequality
$(\ref{Harnack})$ holds for some $\aa>1.$
\item[$(4)$] The log-Harnack inequality $(\ref{Log})$ holds.
\item[$(5)$] For any $\aa>1$, \beq\beg{split}
\label{H1}\int_Mp_t(x,z)\Big(\ff{p_t(x,z)}{p_t(y,z)}\Big)^{\ff 1
{\aa-1}}\mu(\d z) \le
&\exp\Big[\ff{K\aa\rr(x,y)^2}{2(\aa-1)^2(1-\e^{-2Kt})}\Big],\\
&\ \ \ t>0, x,y\in M.\end{split}\end{equation}
\item[$(6)$] There exists $\aa>1$ such that $(\ref{H1})$ holds.
 \item[$(7)$] The following entropy inequality  holds:
\beq\label{H2}\int_Mp_t(x,z)\log\ff{p_t(x,z)}{p_t(y,z)}\mu(\d z)\le
\ff{K\rr(x,y)^2}{2(1-\e^{-2Kt})}, \ \ \ t>0, x,y\in M.\end{equation}\end{enumerate}
\end{thm}

To see that the assumption on the boundary is essential, we intend
to prove that when $\pp M$ is non-empty, each of (\ref{Harnack}),
(\ref{Log}), (\ref{H1}) and (\ref{H2})  implies the convexity of
$\pp M$. Due to technical reasons for estimates on local times, we
assume that $L\rr_\pp$ is bounded for small $\rr_\pp$, where
$\rr_\pp$ is the Riemmanian distance to $\pp M$. This assumption  is
trivial when the manifold is compact. Moreover, by Kasue's
comparison theorems  \cite{K2}, this assumption follows if
there exists $r_0>0$ such that $\<Z,\nn\rr_\pp\>$ is bounded on the set
$\{\rr_\pp\le r_0\}$, $\pp M$ has a bounded second
fundamental form and a strictly positive
 injectivity radius,   the
sectional curvature of $M$ is bounded above,  and the Ricci
curvature of $M$ is bounded below (see e.g. \cite{W05, W07} for
details).

\beg{thm}\label{T1.2} Let  $M$ have a boundary $\pp M$ such that for
some constant $r_0>0$ the function $\rr_\pp$ is smooth with bounded
$L\rr_\pp$ on the set $\{\rr_\pp\le r_0\}$.
 Then $(\ref{Log})$ implies that $\pp M$ is
convex. Consequently, each of statements $(2)$-$(7)$ in Theorem $\ref{T1.1}$
is equivalent to

\ \newline$(8)\ \pp M$ is convex and $(\ref{1.1})$ holds.  \end{thm}

Obviously, Theorem \ref{T1.2} implies the assertions claimed in
Abstract. We remark that a formula for the second fundamental form
was presented in a recent work \cite{W09} for compact manifolds with
boundary by using the gradient estimate due to Hsu \cite{Hsu}. As a
consequence, the manifold is convex if and only if the gradient
estimate

$$|\nn P_t f|^p\le \e^{Kt}P_t |\nn f|^p,\ \ \ t\ge 0, f\in
C_b^1(M)$$ holds for some $p\ge 1$ and $K\in \R.$ When $\pp M$ is
empty  it is well known that such a gradient estimate is equivalent
to the curvature condition (\ref{1.1}) (see e.g. \cite{S}), but the
equivalence with the convexity of boundary  was first observed in
\cite{W09}. Theorem \ref{T1.2} in this paper provides more
equivalent semigroup (heat kernel) properties for (\ref{1.1}) and
the convexity of $\pp M$ without using gradient.

In Section 2 we shall provide in the next section some general
properties for Harnack type inequalities, which are interesting by
themselves. Using these properties we are able to   present complete
proofs for the above two theorems in Sections 3 and 4 respectively.
The log-Harnack inequality is established in Section 5 for  a class
of non-convex manifolds. As an application, the HWI inequality is
presented in Section 6. Finally, two technical points, i.e. the
exponential estimates of the local time and a simple proof of Hsu's
gradient estimate on non-compact manifolds, are addressed in the
Appendix.

\section{Some properties of Harnack Inequalities}
Let $(E,\rr)$ be a metric space, and $P(x,\d y)$ a transition
probability on $E$, which provides a contractive linear operator $P$
on $\scr B_b(E)$, the set of all bounded measurable functions on
$E$:

$$Pf(x)= \int_E f(y)P(x,\d y),\ \ \ f\in \scr B_b(E), x\in E.$$ Let
$\scr B_b^+(E)$ be the set of nonnegative elements in $\scr
B_b(E)$. We shall study the following Harnack inequality with a
power $\aa>1$:

\beq\label{H2'} (Pf(x))^\aa\le (Pf^\aa(y))\exp\Big[\ff{\aa
c\rr(x,y)^2}{\aa-1}\Big],\ \ \ f\in \scr B_b^+(E), x,y\in E,
\end{equation} where $c>0$ is a constant. To state our first result
in this section, we shall assume that $E$ is a length space, i.e.
for any $x\ne y$ and any $s\in (0,1)$, there exists a sequence
$\{z_n\}\subset E$ such that $\rr(x,z_n)\to s\rr(x,y)$ and
$\rr(z_n,y)\to (1-s)\rr(x,y)$ as $n\to\infty.$

\beg{prp}\label{P2.1} Assume that $(E,\rr)$ is a length space and
let $\aa_1,\aa_2>1$ be two constants. If $(\ref{H2'})$ holds for
$\aa=\aa_1,\aa_2,$ it holds also for $\aa=\aa_1\aa_2.$ \end{prp}

\beg{proof} Let

$$s= \ff{\aa_1-1}{\aa_1\aa_2-1},\ \ \ 1-s=
\ff{\aa_1(\aa_2-1)}{\aa_1\aa_2-1},$$  and let $\{z_n\}\subset E$
such that $\rr(x,z_n)\to s\rr(x,y)$ and $\rr(z_n,y)\to
(1-s)\rr(x,y)$ as $n\to\infty.$ Since (\ref{H2'}) holds for
$\aa=\aa_1$ and $\aa=\aa_2,$ for any $f\in \scr B_b^+(E)$ we have

\beg{equation*}\beg{split} &(P f(x))^{\aa_1\aa_2}\le (P
f^{\aa_1}(z_n))^{\aa_2}\exp\Big[\ff{\aa_1\aa_2c\rr(x,z_n)^2}{\aa_1-1}\Big]\\
&\le (P_f^{\aa_1\aa_2}(y))
\exp\Big[\ff{\aa_1\aa_2c\rr(x,z_n)^2}{\aa_1-1}+
\ff{\aa_2c\rr(z_n,y)^2}{\aa_2-1} \Big].\end{split}\end{equation*}
Letting $n\to\infty$ we arrive at

\beg{equation*}\beg{split}(Pf(x))^{\aa_1\aa_2} &\le
(Pf^{\aa_1\aa_2}(y))\exp\Big[\ff{\aa_1\aa_2cs^2\rr(x,y)^2}{\aa_1-1}+
\ff{\aa_2c(1-s)^2\rr(x,y)^2}{\aa_2-1}\Big]\\
&=(Pf^{\aa_1\aa_2}(y))\exp\Big[\ff{\aa_1\aa_2c\rr(x,y)^2}{\aa_1\aa_2-1}\Big].\end{split}\end{equation*}\end{proof}

\beg{prp}\label{P2.2} If $(\ref{H2'})$ holds for some $\aa>1$, then

$$P(\log f)(x)\le \log Pf(y)+ c\rr(x,y)^2,\ \ \ x,y\in E, f\ge 1, f
\in \scr B_b(E).$$\end{prp}

\beg{proof} By Proposition \ref{P2.1}, (\ref{H2}) holds for $\aa^n
(n\in \mathbb N)$ in place of $\aa$. So,

$$P f^{\aa^{-n}}(x)\le (P
f(y))^{\aa^{-n}}\exp\Big[\ff{c\rr(x,y)^2}{\aa^n-1}\Big].$$
Therefore, by the dominated convergence theorem

\beg{equation*}\beg{split} P(\log f)(x)&= \lim_{n\to\infty}
P\Big(\ff{f^{\aa^{-n}}-1}{\aa^{-n}}\Big)(x)\\
&\le \lim_{n\to\infty} \Big\{\ff{(Pf(y))^{\aa^{-n}}-1 }{\aa^{-n}}+
(Pf(y))^{\aa^{-n}}\ff{\exp\big[\ff{c\rr(x,y)^2}{\aa^n-1}\big]-1}{\aa^{-n}} \Big\}\\
&= \log Pf(y)+ c\rr(x,y)^2.\end{split}\end{equation*}\end{proof}

\beg{prp}\label{P2.3} Let $\Phi$ be  a positive function  on
$E\times E$ such that  $\Phi(x,y)\to 0$ as $y\to x$ holds for any
$x\in E.$ Then the log-Harnack inequality

\beq\label{LH} P(\log f)(x)\le \log Pf(y)+\Phi(x,y),\ \ \ x,y\in E,
f\ge 1, f\in \scr B_b(E)\end{equation} implies the strong Feller
property of $P$, i.e. $P\scr B_b(E)\subset C_b(E).$
\end{prp}

\beg{proof} It suffices to prove that $Pf\in C_b(E)$ for $f\in\scr
B_b^+(E).$ Applying (\ref{LH}) for
 $1+\vv f$ in place of $f$,  we obtain

 $$P f(y)-\vv \|f\|_\infty^2 \le P \ff{\log (1+\vv f)}\vv (y)
   \le \ff 1 \vv \log (1+\vv P f(x)) + \ff{\Phi(x,y)}\vv,\ \ \
 \vv>0, x,y\in E.$$
 Letting first $y\to x$  then $\vv\to 0$, we arrive at

 $$\limsup_{y\to x} P f(y) \le P f(x).$$ On the other hand, we have

 $$P\ff{\log (1+\vv f)}\vv (x) -\ff{\Phi(x,y)}\vv
  \le \ff 1 \vv \log (1+\vv P_t f(y))  \le  P_t f(y).$$
  Letting first  $y\to x$ then $\vv\to 0$, we arrive at

  $$Pf(x)\le \liminf_{y\to x} P f(y).$$\end{proof}

  Obviously, each of (\ref{H2'}) and (\ref{LH}) implies that
  $P(x,\cdot)$ and $(P(y,\cdot)$ are equivalent to each other. Indeed, if $P(y,A)=0$ then applying (\ref{H2'}) to $f=1_A$
  or applying (\ref{LH}) to $f=1+n1_A$ and letting $n\to \infty$, we conclude that $P(x,A)=0.$
By the same reason, $P(x,\cdot)$ and $P(y,\cdot)$ are equivalent for any $x, y\in E$  if

  \beq\label{H3} (Pf(x))^\aa\le (Pf^\aa(y))\Psi(x,y),\ \ \ x,y\in E,
  f\in \scr B_b^+(E)\end{equation} holds for some positive function $\Psi$
   on $E\times E$. In these cases let

  $$p_{x,y}(z)= \ff{P(x,\d z)}{P(y,\d z)}$$ be the Radon-Nikodym
  derivative of $P(x,\cdot)$ with respect to $P(y,\cdot).$

  \beg{prp}\label{P2.4} Let $\Phi,\Psi$ be   positive functions on $E\times
  E$.

   $(1)\ (\ref{H3})$ holds if and only if $P(x,\cdot)$ and
   $P(y,\cdot)$ are equivalent and $p_{x,y}$ satisfies

   \beq\label{H4} P\big\{ p_{x,y}^{1/(\aa-1)}\big\}(x)\le
   \Psi(x,y)^{1/(\aa-1)},\ \ \ x,y\in E.\end{equation}

$(2)\ (\ref{LH})$ holds if and only if $P(x,\cdot)$ and
   $P(y,\cdot)$ are equivalent and $p_{x,y}$ satisfies

   \beq\label{H4''} P \{\log p_{x,y}\}(x)\le
   \Phi(x,y),\ \ \ x,y\in E.\end{equation} \end{prp}

   \beg{proof} (1) Applying (\ref{H3}) to $f_n(z):= \{n\land
   p_{x,y}(z)\}^{1/(\aa-1)},\ n\ge 1,$ we obtain

   \beg{equation*}\beg{split} (P f_n(x))^\aa&\le \Psi(x,y) P
   f_n^\aa(y)= \Psi(x,y) \int_E \{n\land
   p_{x,y}(z)\}^{\aa/(\aa-1)}P(y,\d z)\\
   &\le \Psi(x,y)\int_E \{n\land
   p_{x,y}(z)\}^{1/(\aa-1)}P(x,\d z)= \Psi(x,y)
   Pf_n(x).\end{split}\end{equation*}Thus,

   $$P\big\{p_{x,y}^{1/(\aa-1)}\big\}(x)=\lim_{n\to\infty} Pf_n(x)\le
   \Psi(x,y)^{1/(\aa-1)}.$$ So, (\ref{H3}) implies (\ref{H4}).

   On the other hand, if (\ref{H4}) holds then for any $f\in \scr
   B_b^+(E),$ by the H\"older inequality

   \beg{equation*}\beg{split} Pf(x) &= \int_E \{p_{x,y}\}(z) f(z) P(y,\d
   z)\le (Pf^\aa(y))^{1/\aa}\bigg(\int_E p_{x,y}(z)^{\aa/(\aa-1)}
   P(y,\d z)\bigg)^{(\aa-1)/\aa}\\
   &= (Pf^\aa(y))^{1/\aa} (Pp_{x,y}^{1/(\aa-1)}(x))^{(\aa-1)/\aa}\le
   (Pf^\aa(y))^{1/\aa}\Psi(x,y)^{1/\aa}.\end{split}\end{equation*}
   Therefore, (\ref{H3}) holds.

   (2) We shall use the following Young inequality: for any probability
measure $\nu$ on $M$, if $g_1,g_2\ge 0$ with $\nu(g_1)=1$, then

$$\nu(g_1g_2)\le \nu(g_1\log g_1)+\log \nu(\e^{g_2}).$$ For $f\ge
1$, applying the above inequality for $g_1= p_{x,y}, g_2= \log f$
and $\nu= P(y,\cdot)$, we obtain

\beg{equation*}\beg{split} P(\log f)(x) & = \int_E \{p_{x,y}(z)\log
f(z)\} P(y,\d z)\\
&\le P(\log p_{x,y})(x)+ \log Pf(y).\end{split}\end{equation*} So,
(\ref{H4''}) implies (\ref{LH}). On the other hand, applying
(\ref{LH}) to $f_n= 1 + np_{x,y},$ we arrive at

\beg{equation*}\beg{split} &P \{\log p_{x,y}\}(x)\le P(\log f_n)(x) -\log
n\\
&\le \log Pf_n(y)-\log n + \Phi(x,y)= \log \ff{n+1}n
+\Phi(x,y).\end{split}\end{equation*} Therefore, by letting $n\to
\infty$ we obtain (\ref{H4''}).\end{proof}

\section{Proof of Theorem \ref{T1.1}}

By \cite[Lemma 2.2]{W97}, if $\pp M$ is either convex or empty then
 (\ref{1.1}) implies (\ref{Harnack}). Combining
this with Propositions 2.2 and 2.4 for $P=P_t$ so that $p_{x,y}(z)=\ff{p_t(x,z)}{p_t(y,z)}$,
it remains to prove that
(\ref{Log}) implies (\ref{1.1}).

 Let $x\in M$ (when $M$ has a convex
boundary, we take $x$ in the interior) and $X\in T_xM$ be fixed.
 For any $n\ge 1$ we may take $f\in C^\infty_b(M)$ such that $f\ge 1, \ f$ is constant outside a compact set, and

\beq\label{2.2} \nn f(x)=X,\ \ \Hess_f (x)=0,\ \ \ f(x)\ge
n.\end{equation} If $M$ has a convex boundary $\pp M,$ we may assume
further that $f$ is constant in a neighborhood of $\pp M$ so that
the Neumann boundary condition is satisfied.  Such a function can be
constructed by using the exponential map as follows. Let $r_0>0$ be
smaller than the injectivity radius at point $x$ such that the
exponential map

$$\exp_x: \{Y\in T_x M:\ |Y|<r_0\}\to B(x, r_0):= \{z\in M: \rr(x,z)<r_0\}\subset M\setminus \pp M$$
is
diffeomorphic. Then the function

$$g(z):= \<X, \exp_x^{-1}(z)\>,\ \ \ z\in B(x,r_0)$$ is smooth and satisfies
$\nn g(x)=X, \Hess_g(x)=0.$ Let $F\in C_0^\infty(M)$ such that
$F|_{B(x, r_0/4)}=1$ and $F|_{B(x, r_0/2)^c}=0.$ Then $f:= gF +R$
meets our requirements for
 a large enough constant $R>0$.

Taking $\gg_t=\exp_x[-2t\nn\log f(x)],$ we have $\rr(x, \gg_t)= 2 t
|\nn\log f|(x)$ for $t\in [0,t_0],$ where $t_0>0$ is such that
$2t_0|X|<r_0f(x).$ By  (\ref{Log}) with $y=\gg_t$, we obtain

\beq\label{2.3} P_t (\log f)(x) \le \log P_tf(\gg_t)+\ff{2Kt^2|\nn\log
f|^2(x)}{1-\e^{-2Kt}},\ \ \ t\in (0,t_0].\end{equation} Since $Lf\in
C_0^2(M)$ and $L\log f=0$ around $\pp M$, and noting that $\Hess_f(x)=0$ implies $\nn|\nn f|^2(x)=0$, at point $x$ we have

\beg{equation*}\beg{split} \ff{\d}{\d t}P_t\log f|_{t=0}&= L\log f = \ff{Lf}f -|\nn\log f|^2,\\
\ff{\d^2}{\d t^2}P_t\log f|_{t=0} &= L^2\log f= \ff{L^2f}f -\ff{(Lf)^2}{f^2}
+\ff{2|\nn f|^2 Lf}{f^3} + 2 \<\nn Lf, \nn f^{-1}\> -\ff{L|\nn f|^2}{f^2}\\
&\qquad + \ff{2 |\nn f|^2Lf}{f^3}
-\ff{6|\nn f|^4}{f^4} - 2\<\nn |\nn f|^2,\nn f^{-2}\>\\
&=  \ff{L^2f}f -\ff{(Lf)^2}{f^2} -\ff{2}{f^2}\<\nn Lf,\nn f\> -\ff{L|\nn f|^2}{f^2} + \ff{4|\nn f|^2
Lf}{f^3} -\ff{6|\nn f|^4}{f^4}=:A.\end{split}\end{equation*} Thus, by Taylor's expansions,

\beq\label{2.4} P_t(\log f)(x)= \log f(x)+ t\big(f^{-1} Lf -|\nn \log
f|^2\big)(x) +\ff {t^2} 2 A+ o(t^2)\end{equation} holds for small
$t>0.$ On the other hand, let $N_t= //_{x\to \gg_t}\nn\log f(x),$
where $//_{x\to \gg_t}$ is the parallel displacement along the
geodesic $t\mapsto \gg_t$. We have $\dot\gg_t =- 2 N_t$ and
$\nn_{\dot\gg_t} N_t=0.$ So,

\beg{equation*}\beg{split} \ff{\d}{\d t}\log P_t f(\gg_t)|_{t=0}= &\Big(\ff{LP_t f}{P_t f}(\gg_t)
-\ff{2 \<\nn P_t f,N_t\>}{P_t f}(\gg_t)\Big)\Big|_{t=0}= \ff{Lf}f -2|\nn\log f|^2,\\
\ff{\d^2}{\d t^2}\log P_t f(\gg_t)|_{t=0}= &\ff{L^2 f}f -\ff{(Lf)^2}{f^2} -2\<\nn (f^{-1} Lf),
\nn\log f\> -\ff 2 f \<\nn Lf,\nn\log f\> \\
& + \ff 2 {f^2} \<\nn f, \nn\log f\>Lf + 4 \Hess_{\log f}( \nn\log f, \nn\log f)\\
=& \ff{L^2 f}f -\ff{(Lf)^2}{f^2}-4\ff{\<\nn Lf,\nn f\>}{f^2} +
4\ff{|\nn f|^2Lf}{f^3} -4\ff{ |\nn
f|^4}{f^4}=:B,\end{split}\end{equation*} where, as in above, the
functions take value at point $x$ and we have used $\Hess_f(x)=0$ in
the last step. Thus, we have

$$\log P_t f(\gg_t)= \log f(x) + t\big(f^{-1}Lf - 2 |\nn \log f|^2\big)(x) +
\ff {t^2} 2 B +o(t^2).$$ Combining this with (\ref{2.3}) and
(\ref{2.4}), we arrive at

$$\ff 1 t \Big(1- \ff{2Kt}{1-\e^{-2Kt}}\Big)|\nn \log f|^2(x)\le \ff 1 2 \Big(\ff{L|\nn f|^2-
2\<\nn Lf,\nn f\>}{f^2} + \ff{2|\nn f|^4}{f^4}\Big)(x)+o(1).$$
Letting $t\to 0$ we obtain

$$\GG_2(f,f)(x):= \ff 1 2 L|\nn f |^2 (x)-\<\nn Lf, \nn f\>(x)\ge -K|\nn f|^2(x)-
\ff{|\nn f|^4}{f^2}(x).$$
Since  by the Bochner-Weitzenb\"ock formula and (\ref{2.2}) we have $\nn f(x)=X, f(x)\ge n$ and

$$\GG_2(f,f)(x)= \Ric(X,X)-\<\nn_X Z,X\>,$$ it follows that

$$\Ric(X,X)-\<\nn_X Z,X\>\ge -K|X|^2 -\ff{|X|^4} n,\ \ \ n\ge 1.$$ This implies (\ref{1.1})
by letting $n\to\infty.$

\section{Proof of Theorem \ref{T1.2}}

Since in the proofs of \cite[Theorem 2.1 and Lemma 2.2]{W09} only
the boundedness of $L\rr_\pp$  on $\{\rr_\pp\le r_0\}$ rather than
the compactness of $M$ is used, these two results hold true in the
setting of Theorem \ref{T1.2}. More precisely, we have the following
result.

\beg{prp} \label{P3.1} If there exists $r_0>0$ such that $\rr_\pp$
is smooth with bounded $L\rr_\pp$ on $\{\rr_\pp\le r_0\},$ then
there exists a constant $c>0$ such that $\E l_t^2 \le c t$ holds for
  all $x_0\in\pp M$ and $t\in [0,1],$ and

$$ \limsup_{t\to 0} \ff 1 t \Big|\E l_t -\ff 2 {\ss\pi}
\ss{t}\Big|<\infty$$ holds uniformly in $x_0\in\pp M.$\end{prp}

Let $N$ be the unit inward normal vector field of $\pp M$. Then

$$\II(X,X):= -\<\nn_X N, X\>\ge 0,\ \ \ X\in T\pp M$$ is  the second fundamental form of $\pp M$.
By definition $\pp M$ is called convex if $\II\ge 0.$

For any $x\in \pp M$ and $X\in T_x\pp M$, let $f\in C^\infty(M)$ be
such that $f\ge 1, Nf|_{\pp M}=0$ and $\nn f(x)=X.$ We may further
assume that $f$ is constant outside a compact set. To construct such
a function, let $\tt f\in C_0^\infty(\pp M)$ such that $\nn_{\pp
M}\tt f(x)=X,$ where $\nn_{\pp M}$ is the gradient on $\pp M$ with
respect to the induced metric. Let  $\tt f$ be supported on $\pp
M\cap B(x, m)$ for some $m>0$, where $B(x,m)$ is the open geodesic
ball around $x$ with radius $m$. Then there exists $r_1\in (0,1)$
such that
 the exponential map

$$ U:=(B(x,m+3)\cap \pp M)\times [0,r_1)\ni (\theta,r)\mapsto \exp_\theta [rN]$$ is smooth and one-to-one,
which is known as the local polar coordinates around
$B(x,m+2)\cap\pp M$. Let $h\in C^\infty([0,\infty))$ such that
$h|_{[0, (r_1\land r_0)/4]}=1$ and $h|_{[(r_0\land r_1)/2,
\infty)}=0$. Since $\tt f$ is supported on $B(x, m)$  the function

$$M\ni x\mapsto f(x):= R+\beg{cases} \tt f(\theta) h(r), &\text{if\
there \ exists}\ (\theta,r)\in U\ \text{such\ that}\
x=\exp_\theta[rN],\\
0, &\text{otherwise}\end{cases}$$ for large enough constant $R>0$
meets our requirements.

Let $\exp_x^\pp: T_x\pp M\to \pp M$  be the exponential map on the Riemannian manifold $\pp M$ with the induced metric, and let

$$\gg_t= \exp_x^\pp\big[-2t\nn\log f(x)\big],\ \ \ t\ge 0.$$ Applying (\ref{Log})
to $y=\gg_t$ we obtain

\beq\label{F1} P_t \log f(x)\le \log P_t f(\gg_t) + \ff{2K
t^2|\nn\log f|^2(x)}{1-\e^{-2Kt}},\ \ \ t\ge 0.\end{equation} Since
$f$ and $Lf$ satisfy the Neumann boundary condition, we have

\beg{equation}\label{A1}\beg{split} P_t \log f(x)&= \log f(x)
+\int_0^t P_s
L\log f(x)\d s\\
&= \log f(x)+ \int_0^t P_s\ff{Lf} f (x)\d s -\int_0^t P_s|\nn \log
f|^2(x)\d s.\end{split}\end{equation} Let $X_s$ be the reflecting
$L$-diffusion process with $x_0=x$, and let $l_s$ be its local time
on $\pp M$. By the It\^o formula for $|\nn\log f|^2(x_s)$ we obtain

$$P_s |\nn \log f|^2(x)= |\nn \log f|^2(x) +\int_0^sP_r L|\nn\log f|^2(x)\d r +\E\int_0^s \<N, \nn
|\nn\log f|^2\>(X_r)\d l_r.$$
  Since
$f$ satisfies the Neumann boundary condition so that

$$\<N, \nn|\nn \log f|^2\> = 2f^{-2} \Hess_f(N, \nn f),$$
and since $\<\nn f,\nn \<N, \nn f\>\>=0$ implies

$$\Hess_f(N,\nn f)= - \<\nn_{\nn f}N,\nn f\>= \II(\nn f,\nn f),$$ it
follows that

$$P_s|\nn\log f|^2(x)= |\nn\log f|^2(x) + O(s)+ 2f^{-2}(x)\II(\nn f,\nn
f)(x) \E l_s + o(\E l_s).$$ Since due to Proposition \ref{P3.1} we
have $\lim_{t\to 0} t^{-1/2} \E l_t = \ff 2 {\ss\pi},$ this and
(\ref{A1}) yield (recall that $\nn f(x)=X$)

\beq\label{F2}P_t\log f(x)= \log f(x)+ \int_0^t P_s\ff{Lf}f(x)\d s
-|\nn \log f|^2(x) - \ff {8 t^{3/2}}{3\ss\pi\, f^2(x)}\II(X,X)
+o(t^{3/2}).\end{equation} On the other hand, we have

\beg{equation*}\beg{split} P_t f(\gg_t)&= f(\gg_t) + \int_0^t P_s
Lf(\gg_t)\d s\\
&= f(x)+t\<\dot \gg_s,\nn f(\gg_s)\>|_{s=0} + O(t^2) +\int_0^t P_s
Lf(x) \d s\\
&= f(x)-\ff{2t}{f(x)} |\nn f|^2(x) +\int_0^t P_s Lf(x)\d s
+O(t^2).\end{split}\end{equation*} Thus,

$$\log P_t f(\gg_t)= \log f(x) + \ff 1 {f(x)}\int_0^t P_s Lf(x)\d s- 2t|\nn\log f|^2(x) + O(t^2).$$
Combining this with (\ref{F1}) and (\ref{F2}) we arrive at

\beg{equation}\label{F3}\beg{split} & \ff 1 {t\ss t } \int_0^t
\Big(P_s\ff{Lf}f - \ff{P_sLf} f\Big)(x)\d s +\ff 1 {\ss t}\Big(
1-\ff{2Kt}{1-\e^{-2Kt}}\Big)|\nn \log f|^2(x)\\
&\le \ff 8{3\ss\pi\, f^2(x)} \II(X,X) +o(1
).\end{split}\end{equation} Obviously,

$$\lim_{t\to 0}\ff 1 {\ss t} \Big(1-\ff{2Kt}{1-\e^{-2Kt}}\Big)=0.$$ So, to derive
$\II(X,X)\ge 0$ from (\ref{F3}) it remains to verify

\beq\label{F4} \lim_{t\to 0}\ff 1 {t\ss t } \int_0^t
\Big(P_s\ff{Lf}f - \ff{P_sLf} f\Big)(x)\d s=0.\end{equation} Noting
that $Z$ is $C^2$-smooth and $f\in C^\infty(M)$ is constant outside
a compact set, we have $Lf\in C_0^2(M)$. Moreover, $f\ge 1$ and
$f$ satisfies the Neumann boundary condition. So, by the It\^o formula
we have

\beq\label{F5} \beg{split} &\Big(P_s\ff {Lf}f-
\ff{P_sLf}f\Big)(x)\\
&= \int_0^s \Big(P_r L\ff {Lf}f - \ff{P_r L^2 f} f\Big)(x)\d r +\E
\int_0^s\Big(\ff 1 {f(X_r)}- \ff 1 {f(x)}\Big)\<N, \nn Lf\>(X_r) \d
l_r.\end{split}\end{equation} Since $\ff 1 f$ is bounded and $X_r\to
x$ as $r\to 0$, it follows from Proposition \ref{P3.1} that

\beg{equation*}\beg{split} &\limsup_{s\to 0}\ff 1{\ss s}
\bigg|\E\int_0^s\Big(\ff 1{f(x_r)}- \ff 1 {f(x)}\Big)\<N, \nn
Lf\>(x_r) \d l_r\bigg|\\
&\le \limsup_{s\to 0} \ff{\|\nn Lf\|_\infty}{\ss s}
\E\Big(l_s\sup_{r\in [0,s]} |f(X_r)^{-1}- f(x)^{-1}|\Big)\\
&\le \|\nn Lf\|_\infty \limsup_{s\to 0} \Big(\ff{\E
l_s^2}s\Big)^{1/2} \Big(\E\sup_{r\in [0,s]} |f(X_r)^{-1}-
f(x)^{-1}|^2\Big)^{1/2}=0.\end{split}\end{equation*} Therefore, (\ref{F4})
follows from (\ref{F5}) immediately.

\section{An extension to non-convex manifolds}

In this section we aim to established  the log-Harnack inequality on a class of non-convex manifolds.
To this end, we need the following assumption.

\paragraph{(A)}  \emph{The boundary $\pp M$ has a bounded second fundamental form and a strictly positive
 injectivity radius,   the
sectional curvature of $M$ is bounded above,    and there exists
$r>0$ such that $Z$ is bounded on the $r$-neighborhood of $\pp M$.
}

Under this assumption, we have $\sup_x\in M \e^{\ll l_t}<\infty$ for
all $\ll>0$ (see Proposition \ref{AA1} in Appendix). Let

$$U_{x,y}(s)= \sup_{z: \rr(z,x)\lor \rr(z,y)\le \rr(x,y)} \mathbb E^z \e^{2\si l_s},\ \ \ x,y\in M, s\ge 0.$$

As a complement to known equivalent statements
for lower bounds on  curvature and second fundamental form derived recently in \cite{W09b},
  the following result  provides two more
equivalent statements.

\beg{thm} Assume {\bf (A)}. Let $K,\si\in \R$ be two constants.
Then the following statements are
equivalent each other:
\begin{enumerate}
\item[$(1)$] \ $\Ric-\nn Z\ge -K, \II\ge -\si.$
\item[$(2)$]\ $P_t (\log f)(x)\le \log P_t f(y)
+\ff{\rr(x,y)^2}{4\int_0^t\e^{-2Ks}\big\{U_{x,y}(s)\big\}^{-1}\d s}$
holds for all   $ f\in \B_b^+(M)$ with $f\ge 1, t\ge 0,$ and $x,y\in
M.$
\item[$(3)$]\ $ \int_M p_t(x,z)\log \ff{p_t(x,z)}{p_t(y,z)}\mu(\d
z)\le\ff{\rr(x,y)^2}{4\int_0^t\e^{-2Ks}\big\{U_{x,y}(s)\big\}^{-1}\d
s}$ holds for all $t>0, x,y\in M.$\end{enumerate}\end{thm}

\beg{proof} Since Proposition \ref{P2.4} ensures that (2) and (3) are equivalent, it suffices to prove the equivalence of (1) and (2).

(a) (1) implies (2). According to  (1), the following Hsu's gradient estimate holds (see Proposition \ref{A2} in Appendix):

\beq\label{Fg1} |\nn P_t f|^2 \le \big(\E\big\{|\nn f|(X_t)\e^{Kt+\si l_t}\big\}\big)^2\le (P_t|\nn f|^2)\E\e^{2Kt+2\si l_t}.
\end{equation}

Let $\gg: [0,1]\to M$ be the minimal curve with constant such that $\gg(0)=y$ and $\gg(1)=x.$ We have
$|\dot \gg|=\rr(x,y).$ Let $h\in C^1([0,t])$ be such that $h(0)=0$ and $h(t)=1.$ By (\ref{Fg1}) and the definition of $U_{x,y}$  we have

\beg{equation*}\beg{split} &\ff{\d}{\d s} P_s \log P_{t-s}f(\gg\circ h(s))\\
&= -P_s|\nn\log P_{t-s}f|^2(\gg\circ h(s))+\dot h(s)\<\dot \gg\circ h(s),\nn P_s\log P_{t-s}f(\gg\circ h(s))\>\\
&\le -P_s|\nn\log P_{t-s}f|^2(\gg\circ h(s)) + |\dot h(s)|\rr(x,y) \e^{-Ks}\big\{U_{x,y}(s)P_s |\nn \log P_{t-s}f|^2(\gg\circ h(s))\big\}^{1/2}\\
&\le \ff 1 4 |\dot h(s)|^2\rr(x,y)^2 U_{x,y}(s)\e^{2Ks},\ \ \ \ s\in [0,t].\end{split}\end{equation*} This implies

$$P_t\log f(x)\le \log P_t f(y) +\ff {\rr(x,y)^2} 4\int_0^t |\dot h(s)|^2 U_{x,y}(s) \e^{2Ks}\d s.$$ Therefore,
we prove (2) by taking

$$h(s)= \ff{\int_0^s \e^{-2Kr}\{U_{x,y}(r)\}^{-1}\d r}{\int_0^t \e^{-2Kr}\{U_{x,y}(r)\}^{-1}\d r},\ \ \ \ s\in [0,t].$$

(b) (2) implies (1).  Let $x\in M\setminus \pp M.$ There exists $\dd>0$ such that the closed geodesic ball  $\bar B(x, 2\dd)$ at $x$ with
radius $2\dd$ is contained in $M\setminus \pp M,$ i.e. $\bar B(x,2\dd)\cap \pp M=\emptyset.$ Let $\tau$ be the hitting time
of $X_t$ to the boundary, we have (cf. \cite[Proposition A.2]{W09})

$$\P^z(\tau\le t)\le C\e^{-\dd^2/(16t)},\ \ \ z\in B(x,\dd)$$  for some constant $C>0$ and all $t>0.$
Moreover, by \cite[Proof of Lemma 2.1]{W05}, we have

\beq\label{Fg2} C' :=\sup_{z\in\pp M} \E^z \e^{2\si l_1}<\infty.\end{equation} Since $l_t=0$ for $t\le \tau$ and $l_t$ is increasing in $t$, it follows that

\beg{equation*}\beg{split} \E\e^{2\si l_t} &\le \P(\tau>t)+\E1_{\{\tau\le t\}}\E^{X_{\tau}}\e^{2\si l_t}\\
&\le 1 +C' C\e^{-\dd^2/(16t)},\ \ \ \ t\in [0,1].\end{split}\end{equation*} Thus, for any $y\in B(x,\dd),$

$$\int_0^t \ff{\e^{-2Ks}}{U_{x,y}(s)}\d s=\int_0^t \e^{-2Ks}\d s + o(t^3)= \ff{\e^{2Kt}-1}{2K}+o(t^3),$$ where $o(t^3)$ is
uniform in $y\in B(x,\dd).$
Combining this with the proof of Theorem \ref{T1.1}, we derive $\Ric-\nn Z\ge - K$ from (2).

Now, let $x\in \pp M.$ By Proposition \ref{P3.1} and (\ref{Fg2}) we have

$$\sup_{z\in M} \E^z \e^{2\si l_t} \le 1 +\ff{4\si}{\ss\pi} \ss t +O(t).$$ Then

$$\int_0^t \e^{2Ks} \{U_{x,y}(s)\}^{-1}\d s \ge t +\ff {4\si}{\ss\pi}\int_0^t \ss s \d s + o(t^{3/2})=t+ \ff{8\si}{3\ss\pi}
t^{3/2} + o(t^{3/2}).$$ So, (2) implies

$$P_t \log f(x)\le \log P_t f(y) +\ff{\rr(x,y)^2}{1+\ff{8\si}{3\ss\pi} t^{3/2} + o(t^{3/2})}.$$ Thus,
instead of (\ref{F3}) the proof of Theorem \ref{T1.2} yields

\beg{equation*}\beg{split} & \ff 1 {t\ss t } \int_0^t
\Big(P_s\ff{Lf}f - \ff{P_sLf} f\Big)(x)\d s +\ff 1 {t\ss t}\Big(t-t-\ff{8\si}{3\ss\pi}t^{3/2} +o(t^{3/2})\Big)|\nn \log f|^2(x)\\
&\le \ff 8{3\ss\pi\, f^2(x)} \II(X,X) +o(1).\end{split}\end{equation*} By this and (\ref{F4}) and letting $t\to 0$ we deduce
that $\II(X,X)\ge -\si|X|^2.$
\end{proof}

\section{HWI inequality}

To study the HWI inequality,  we consider the symmetric case that
$Z=\nn V$ for some $V\in C^2(M)$ such that $\mu(\d x)=\ \e^{V(x)}\d
x$ is a probability measure on $M$, where $\d x$ is the Riemannian
volume measure on $M$. Let $P_t$ be the semigroup of the reflecting
diffusion process generated by $L$ on $M$, which is  then symmetric
in $L^2(\mu)$. When $\pp M$ is convex (\ref{1.1}) implies the
following gradient estimate (cf. \cite{Q, W97})

\beq\label{G} |\nn P_t f|\le \e^{Kt}P_t|\nn f|,\ \ \ f\in
C_b^1(M).\end{equation} Combining this estimate and an argument of
\cite{BGL} (see also \cite{OV}), we can easily obtain the following HWI inequality:

\beq\label{HWI} \mu(f^2\log f^2)\le 2 \ss{\mu(|\nn f|^2)}\,
W_2(f^2\mu,\mu)+\ff K 2 W_2(f^2\mu, \mu)^2,\ \ \
\mu(f^2)=1,\end{equation} where $W_2$ is the $L^2$-Wasserstein
distance induced by the Riemannian distance function $\rr$ on $M$.
More precisely, for a probability measure $\nu$ on $M$  (note that
we are using $\rr^2$ to replace $\ff 1 2 \rr^2$ in \cite{BGL})

$$W_2(\nu,\mu)^2:= \inf_{\pi\in\scr C(\nu,\mu)} \int_{M\times M} \rr(x,y)^2\pi(\d x,\d y),$$ where
$\scr C(\nu,\mu)$ is the class of all couplings of $\nu$ and $\mu.$

\beg{thm}\label{T6.2} Let $Z=\nn V$ for some $V\in C^2(M)$ such that
$\mu$ is a probability measure. Assume {\bf (A)} and $(\ref{1.1})$.
Let $\II\ge -\si$ for some $\si\in \R$. Then

$$\eta_\ll(s):= \sup_{x\in M}\mathbb E^x \e^{\ll l_s}<\infty,\ \ \ s,\ll\ge 0$$ holds, and for any $t>0$,

\beq\label{HWI2} \mu(f^2\log f^2)\le 4\bigg(\int_0^t \e^{2Ks}\eta_{2\si}(s)\d s\bigg)\mu(|\nn f|^2)
+ \ff {W_2(f^2\mu,\mu)^2} {4\int_0^t\e^{-2Ks}\eta_{2\si }(s)^{-1} \d s},\ \ \mu(f^2)=1.\end{equation}
\end{thm}

\beg{proof} By  Proposition  \ref{AA1} in Appendix, it remains to
verify (\ref{HWI2}). Let $f\in C_b^1(M)$ and $t>0$. We have

\beq\label{3.3} \ff{\d}{\d s} P_s\big\{(P_{t-s}f^2)\log
P_{t-s}f^2\big\}=P_s\ff{|\nn P_{t-s}f^2|^2}{P_{t-s}f^2},\ \ \ s\in
[0,t].\end{equation} By Proposition  \ref{A2} below and the Schwartz inequality we
have

\beg{equation*}\beg{split} \ff{|\nn
P_{t-s}f^2|^2}{P_{t-s}f^2}(y)&\le \e^{2K(t-s)}\ff{(\EE^y\{|\nn
f^2|(X_{t-s})\e^{\si l_{t-s}}\})^2}{P_{t-s} f^2(y)}\\
&\le 4 \e^{2K(t-s)} \EE^y \{|\nn f|^2(X_{t-s})\e^{2\si
l_{t-s}}\}=:4\e^{2K(t-s)}g_s(y),\ \ \ s\in [0,t], y\in
M.\end{split}\end{equation*} Combining this with (\ref{3.3}) we
obtain

$$P_t (f^2\log f^2)\le (P_tf^2)\log P_t f^2 + 4 \int_0^t
\e^{2K(t-s)}P_s g_s\d s.$$Since $\mu$ is an invariant measure of
$P_t$, taking integral for both sides with respect to $\mu$ we
arrive at

\beq\label{AAA}\mu(f^2\log f^2) \le \mu((P_tf^2)\log P_t f^2)+ 4\int_0^t
\e^{2K(t-s)}\mu(g_s)\d s.\end{equation} Let $P_t^\si$ be defined by

$$P_t^\si h(x)= \mathbb E^x[h(X_t)\e^{2\si l_t}],\ \ h\in C_b(M).$$ Then it is easy to see that
$u(t,x):= P_t^\si h(x)$ solve the heat equation with Robin boundary condition

$$\pp_t u= Lu,\ \ u(0,\cdot)=h, (Nu+2\si u)|_{\pp M}=0.$$ In particular, since $L$ is symmetric in $L^2(\mu)$ under
the Robin boundary  condition, so is $P_t^\si$. Therefore,

$$\mu(g_s)= \\mu(P_{t-s}^\si |\nn f|^2)= \mu (|\nn f|^2 P_{t-s}^\si 1) \le
\mu(|\nn f|^2) \eta_{2\si}(t-s).$$  Combining this with (\ref{AAA}) we obtain

\beq\label{3.4} \mu(f^2\log f^2) \le \mu((P_tf^2)\log P_t f^2)+
4\mu(|\nn f|^2 ) \int_0^t \e^{2Ks}\eta_{2\si}(s)\d s.\end{equation}

On the other hand, for any $x,y\in M$, let $x_\cdot:
[0,1]\to M$ be the minimal curve linking $x$ and $y$ with constant
speed. We have $|\dot x_s|=\rr(x,y)$. Let $h\in C^1([0,t])$ be such
that $h_0=1, h_t=0.$ Then by Proposition \ref{A2} below, we have

\beg{equation}\label{3.5}\beg{split} &P_t\log f^2(x)-\log P_t
f^2(y)= \int_0^t
\ff{\d }{\d s} P_s(\log P_{t-s}f^2)(x_{h_{t-s}})\d s\\
&\le \int_0^t \Big\{|\dot h_{t-s}| \rr(x,y) |\nn P_s(\log
P_{t-s}f^2)|(x_{h_{t-s}})-\EE^{x_{h_{t-s}}} \ff{|\nn
P_{t-s}f^2|^2}{(P_{t-s}f^2)^2}(X_s)  \Big\}\d s\\
&\le \int_0^t \EE^{x_{h_{t-s}}} \Big\{ |\dot h_{t-s}| \rr(x,y)
\ff{|\nn P_{t-s}f^2|}{P_{t-s}f^2}(X_s)\e^{K(t-s)+\si
l_{t-s}}-\ff{|\nn
P_{t-s}f^2|^2}{(P_{t-s}f^2)^2}(X_s)\Big\}\d s\\
&\le \ff {\rr(x,y)^2}4 \int_0^t {\dot h}_s^2
\e^{2Ks}\eta_{2\si}(s)\d s=:c(t)\rr(x,y)^2.\end{split}\end{equation}

Now, let $\mu(f^2)=1$ and $\pi\in \scr C(f^2\mu,\mu)$ be the optimal
coupling for $W_2(f^2\mu,\mu).$ It follows from the symmetry of
$P_t$ and (\ref{3.5}) that

\beg{equation*}\beg{split} \mu((P_tf^2)\log P_tf^2)&=\mu(f^2 P_t\log
P_t f^2)=\int_{M\times M} P_t(\log P_t f^2)(x)\pi(\d x,\d y)\\
&\le \int_{M\times M} \big\{\log
P_{2t}f^2(y)+c(t)\rr(x,y)^2\big\}\pi(\d x,\d y)\\
&=\mu(\log P_{2t}f^2)+c(t)W_2(f^2\mu,\mu)^2\le
c(t)W_2(f^2\mu,\mu)^2,\end{split}\end{equation*} where in the last
step we have used the Jensen inequality that

$$\mu(\log P_{2t}f^2)\le \log \mu(P_{2t} f^2)=0.$$ Combining this with
(\ref{3.4}) we obtain

$$\mu(f^2\log f^2) \le 4\mu(|\nn f|^2) \int_0^t
\e^{2Ks}\eta_{2\si}(s)\d s+ \ff{W_2(f^2\mu,\mu)^2} 4 \int_0^t {\dot
h}_s^2 \e^{2Ks} \eta_{2\si}(s)\d s.$$ Then the  proof is completed
by taking

$$h_s= \ff{\int_s^t \e^{-2Ku}\eta_{2\si}(u)^{-1}\d u}{\int_0^t
\e^{-2Ku}\eta_{2\si}(u)^{-1}\d u},\ \ \ s\in [0,t].$$ \end{proof}

\section{Appendix}

We aim to confirm the exponential integrability of the local time and Hsu's gradient estimate used in Section 5 and Section 6 for the non-convex case,
which are known in \cite{W05} and \cite{Hsu} respectively for the compact case. Here we shall reprove them for the   non-compact case under
assumption {\bf (A)}.

 To estimate $\E \e^{\ll l_t}$ for $\ll>0$, we introduce some concrete conditions in terms of assumption {\bf (A)}.
 Let $\Sect_M$ be the
sectional curvature of $M$ and $\i_{\pp M}>0$ be the injectivity radius of $\pp M$. Let
$$\dd_r(Z):= \sup_{\pp_r M}\<Z, \nn\rr_{\pp M}\>^-,\ \ \ r>0.$$

\beg{prp} \label{AA1}   Let $r_0,\si, k,
>0$ be such that $\dd_{r_0}(Z)<\infty, -\si\le\II\le \gg$ and
$\Sect_M\le k$.
Then

  $$\sup_{x\in M} \EE^x \e^{\ll l_t} \le \exp\Big[ \ff{\ll d r} 2
  +\Big(\ff{ \ll d} r + \ll \dd_r (Z)+ 2 \ll^2\Big)t\Big],\ \ \
  t\ge 0, \ll\ge 0$$ holds for any

 $$0<r\le \min\bigg\{\i_{\pp M},\ r_0,\ \ff1 {\ss k} \arcsin\bigg(\ff{\ss k}{\ss{k+\gg^2}}\bigg)\bigg\}.$$ \end{prp}

 \beg{proof} Let

 $$h(s)= \cos\big(\ss k\, s\big)- \ff{\gg } {\ss k}\sin\big(\ss k\, s\big),\
 \ \ s\ge 0.$$ Then $h$ is the unique solution to the equation

 $$h'' + k h=0,\ \ \ h(0)=1, h'(0)= -\gg.$$
 By the Laplacian comparison theorem for $\rr_{\pp M}$ (cf.
 \cite[Theorem 0.3]{Kasue} or \cite{W07}),

 $$ \DD \rr_{\pp M} \ge \ff{(d-1)h'}{h}(\rr_{\pp
 M}),\ \ \ \rr_{\pp M}< \i_{\pp M} \land h^{(-1)}(0).$$
 Thus,

 \beq\label{L0} L\rr_{\pp M}\ge \ff{(d-1)h'}{h}(\rr_{\pp M})
 -\dd_r(Z),\ \ \ \rr_{\pp M}\le r.\end{equation} Now, let

 \beg{equation*}\beg{split} &\aa= (1-h(r))^{1-d}\int_0^r
 (h(s)-h(r))^{d-1}\d s,\\
 &\psi(s) =\ff 1 \aa\int_0^s (h(t)-h(r))^{1-d}\d t\int_{t\land
 r}^r(h(u)-h(r))^{d-1}\d u,\ \ s\ge 0.\end{split}\end{equation*} We
 have $\psi(0)=0, 0\le \psi'\le \psi'(0)=1.$ Moreover, as observed in \cite[Proof of
 Theorem 1.1]{W05},

 \beq\label{L1} \aa\ge \ff r d,\ \ \psi(\infty)=\psi(r)\le
 \ff{r^2}{2\aa}\le \ff{dr} 2.\end{equation} Combining this with
 (\ref{L0}) we obtain (note that $\psi'(s)=0$ for $s\ge r$)

 \beq\label{L2} L\psi\circ\rr_{\pp M} =\psi'\circ\rr_{\pp
 M}L\rr_{\pp M}
 +\psi''\circ\rr_{\pp M}\ge -\ff 1 \aa -\dd_r(Z)\ge -\ff d r
 -\dd_r(Z).  \end{equation} On the other hand, since $\psi'(0)=1$, by the It\^o
 formula we have

 \beq\label{CC}\d \psi\circ\rr_{\pp M}(X_t) =\ss 2 \psi'\circ\rr_{\pp M}(X_t) \d
 b_t +L\psi\circ\rr_{\pp M}(X_t)\d t +\d l_t,\end{equation}
 where $b_t$ is the one-dimensional Brownian motion. Then it follows from
 (\ref{L1}) and (\ref{L2}) that (note that $|\psi'|\le 1$)

\beg{equation*}\beg{split} \EE\e^{\ll l_t} &=\EE \exp\bigg[\ll
\psi\circ\rr_{\pp M}(X_t) +\Big(\ff{d\ll} r+\ll \dd_r(Z)\Big)t-\ss 2
\ll \int_0^t\psi'\circ\rr_{\pp M}(X_s)\d b_s\bigg]\\
&\le \exp\Big[\ff 1 2 \ll d r +\Big(\ff{d\ll} r
+\ll\dd_r(Z)\Big)t\Big]
\bigg(\EE\exp\bigg[4\ll^2\int_0^t\big(\psi'\circ\rr_{\pp
M}(X_s)\big)^2\d s\bigg]\bigg)^{1/2}\\
&\le \exp\bigg[\ff 1 2 \ll d r +\Big(\ff {d\ll} r +\ll \dd_r(Z) + 2
\ll^2\Big)t\bigg].\end{split}\end{equation*}
 \end{proof}

\beg{prp}\label{A2} Assume that {\bf (A)}.  Let $\kk_1, \kk_2\in
C_b(M)$ be such that

\beq\label{2.0} \Ric -\nn Z\ge -\kk_1,\ \ \ \II\ge
-\kk_2\end{equation}  hold on $M$ and $\pp M$ respectively. Then

\beq\label{Hsu} |\nn P_t f|(x) \le \EE^x \bigg\{|\nn f|(X_t)
\exp\bigg[ \int_0^t \kk_1(X_s)\d s+\int_0^t\kk_2(X_s)\d
l_s\bigg]\bigg\}\end{equation}
 holds for all $f\in C_b^1(M), t>0, x\in M.$ \end{prp}

We first provide a simple proof of (\ref{Hsu}) under a further  condition that
 $|\nn P_\cdot f|$ is bounded on $[0,T]\times M$ for any $T>0$, then drop this assumption by an approximation argument.
 Since this condition is trivial for compact $M$, our proof below is much
 shorter than that in \cite{Hsu}.

 \beg{lem} \label{L2.1} Assume that $f\in C_b^1(M)$ such that $|\nn P_\cdot f|$ is bounded on
 $[0,T]\times M$ for any $T>0$. Then $(\ref{Hsu})$ holds. \end{lem}

 \beg{proof}  For any $\vv>0$, let

 $$\zeta_s=\ss{\vv +|\nn P_{t-s}f|^2}\, (X_s),\ \ \ s\le t.$$ By the
 It\^o formula we have

\beg{equation*}\beg{split} \d \zeta_s = & \d M_s+\ff{L|\nn
P_{t-s}f|^2 -2 \<\nn LP_{t-s}f, \nn P_{t-s}f\>}
{2\ss{\vv+|\nn P_{t-s}f|^2)^2}}(X_s) \d s\\
&- \ff{|\nn|\nn P_{t-s}f|^2|^2}{4(\vv+|\nn
P_{t-s}f|^2)^{3/2}}(X_s)\d s+\ff{N |\nn P_{t-s}f|^2}{2\ss{\vv+|\nn
P_{t-s}f|^2}} (X_s)\d l_s, \ \ \ s\le t,
\end{split}\end{equation*} where $M_s$ is a local martingale.  Combining this with (\ref{2.0}) and
(see \cite[(1.14)]{Ledoux})
\beq\label{Ledoux}L|\nn u|^2-2\<\nn Lu, \nn u\>\ge -2 \kk_1|\nn u|^2
+\ff{|\nn|\nn u|^2|^2}{2|\nn u|^2},\end{equation}   we obtain

$$\d\zeta_s\ge \d M_s -\ff{\kk_1 |\nn
P_{t-s}f|^2}{\vv+|\nn P_{t-s} f|^2}(X_s)\zeta_s\d s -\ff{\kk_2  |\nn
P_{t-s}f|^2}{\vv + |\nn P_{t-s}f|^2}(X_s)\zeta_s\d l_s,\ \ \ s\le
t.$$ Since   $\zeta_s$ is bounded on $[0,t],$ $\kk_1$ and $\kk_2$
are bounded, and by Proposition \ref{AA1} below $\EE \e^{\ll
l_t}<\infty$ for all $\ll>0$, this implies that

$$[0,t]\ni s\mapsto \zeta_s\exp\bigg[\int_0^s\ff{\kk_1|\nn P_{t-r}f|^2}{\vv+|\nn P_{t-r}f|^2}(X_r)\d r +\int_0^s
\ff{\kk_2|\nn P_{t-r}f|^2}{\vv+ |\nn P_{t-r} f|^2}(X_r)\d
l_r\bigg]$$ is a submartingale for any $\vv>0$. Letting
$\vv\downarrow 0$ we conclude that

$$[0,t]\ni s\mapsto |\nn P_{t-s}f|(X_s) \exp\bigg[\int_0^s\kk_1(X_r)\d r +\int_0^s\kk_2(X_r)\d
l_r\bigg]$$ is a submartingale as well. This completes the
proof.\end{proof}

 By Lemma \ref{L2.1}, to prove Proposition \ref{A2}  it suffices to confirm the boundedness
of $|\nn P_\cdot f|$ on $[0,T]\times M$ for $f\in C_b^1(M).$   Below
we first consider $f\in C_0^\infty(M)$ satisfying the Neumann
boundary condition.

\beg{lem}\label{L2.3} Assume {\bf (A)}. If $(\ref{1.1})$ holds then
for any $T>0$ and $f\in C_0^\infty(M)$ such that $Nf|_{\pp M}=0$,
$|\nn P_\cdot f|$ is bounded on $[0,T]\times M.$\end{lem}

\beg{proof} We shall take a conformal change of metric as in
\cite{W07} to make the boundary convex, so that the known  estimates
 for the convex case can be applied. As explained  on page 1436 in
 \cite{W07}, under assumption {\bf (A)} there exists $\phi\in
 C^\infty(M)$ and a constant $R>1$ such that $1\le\phi\le R, |\nn
 \phi|\le R, N\log\phi |_{\pp M}\ge\si,$ and $\nn\phi=0$ outside
 $\pp_r M.$ Since $\II\ge -\si$, by \cite[Lemma 2.1]{W07} $\pp M$ is
 convex under the new metric

 $$\<\cdot,\cdot\>=\phi^{-2}\<\cdot,\cdot\>.$$ Let $\DD', \nn',\Ric'$
 be corresponding to the new metric. By \cite[Lemma 2.2]{W07}

 $$L' :=\phi^2 L= \DD' + (d-2)\phi\nn\phi +\phi^2 Z=: \DD' +Z'.$$
Following e.g.  \cite{W07} we shall now calculate the curvature
tensor $\Ric'-\nn' Z'$ under
 the new metric. By \cite[(9)]{W07}, for any unit vector $U\in TM$, $U':= \phi  U$
is unit under the new metric, and the corresponding Ricci curvature
satisfies

\beq\label{DD}\beg{split}\Ric'(U',U')\ge &\phi^2 \Ric(U,U)
+\phi\DD\phi -(d-3) |\nn
\phi|^2\\
&-2(U\phi)^2 + (d-2)\phi\Hess_\phi(U,U).\end{split}\end{equation}
 Noting that

$$\nn_X'Y= \nn_X Y- \<X, \nn \log \phi\>Y -\<Y,\nn \log\phi\>X +\<X,Y\>\nn\log \phi,\ \ \ X,Y\in TM,$$
we have

\beg{equation*} \beg{split}&\<\nn_{U'}Z', U'\>' = \<\nn_U Z', U\> -\<Z',\nn\log \phi\>\\
&=\phi^2 \<\nn_U Z,U\> +(U\phi^2) \<Z,U\> +(d-2)(U\phi)^2+(d-2)\phi
\Hess_{\phi} (U,U)-\<Z',\nn\log\phi\>.
\end{split}\end{equation*} Combining this with (\ref{DD}), (\ref{1.1}), $\|Z\|_r<\infty$
and the properties of $\phi$ mentioned above,
 we find a constant $K'\ge 0$ such that

$$\Ric'(U,'U')-\<\nn'_{U'}Z' ,U'\>'\ge -K',\ \ \ \<U',U'\>'=1.$$For
any $x,y\in M$, let $(X_t', Y_t')$ be the coupling by parallel
displacement of the reflecting diffusion processes generated by $L'$
with $(X_0',Y_0')=(x,y).$ Let $\rr'$ be the Riemannian distance
induced by $\<\cdot,\cdot\>'$. Since $(M, \<\cdot,\cdot\>')$ is
convex, we have (see \cite[(3.2)]{W97})

$$\rr'(X_t', Y_t')\le \e^{K't}\rr'(x,y),\ \ \ t\ge 0.$$ Since $1\le
\phi\le R$, we have $R^{-1}\rr\le \rr'\le\rr$ so that

\beq\label{L4} \rr(X_t',Y_t')\le R \e^{K't} \rr(x,y),\ \ \ t\ge
0.\end{equation} To derive the gradient estimate of $P_t$, we shall
make time changes

$$\xi_x(t)= \int_0^t \phi^2(X_s')\d s,\ \ \ \xi_y(t)=
\int_0^t\phi^2(Y_s')\d s.$$ Since $L'=\phi^2 L,$ we see that $X_t:=
X_{\xi_x^{-1}(t)}'$ and $Y_t:= Y'_{\xi_y^{-1}(t)}$ are generated by
$L$ with reflecting boundary. Again by $1\le \phi\le R$ we have

$$R^{-2} t\le \xi_x^{-1}(t), \xi_y^{-1}(t)\le t,\ \ \ t\ge 0.$$
Combining this with $|\nn \phi|\le R, 1\le \phi\le R$ and (\ref{L4})
we arrive at

\beg{equation}\label{L5}\beg{split}
|\xi_x^{-1}(t)-\xi_y^{-1}(t)|&\le \int_{\xi_x^{-1}(t)\land
\xi_y^{-1}(t)}^{\xi_x^{-1}(t)\lor \xi_y^{-1}(t)} \phi^2(Y_s')\d s=
|\xi_y\circ\xi_y^{-1}(t)-\xi_y\circ\xi_x^{-1}(t)|\\
&=|\xi_x\circ\xi_x^{-1}(t)-\xi_y\circ\xi_x^{-1}(t)|\le
\int_0^{\xi_x^{-1}(t)} |\phi^2(X_s')-\phi^2(Y_s')|\d s\\
&\le 2R^2 \rr(x,y) \int_0^t\e^{K's}\d s \le 2t\e^{K't} R^2
\rr(x,y).\end{split}\end{equation} Therefore,

\beq\label{L6} \beg{split} &|P_tf(x)-P_tf(y)| =
|\EE\{f(X_{\xi_x^{-1}(t)}')- f(Y'_{\xi_y^{-1}(t)})\}|\\
&\le \EE |f(X'_{\xi_y^{-1}(t)})-f(Y'_{\xi_y^{-1}(t)})| +
|\EE\{f(X'_{\xi_x^{-1}(t)})-f(X'_{\xi_y^{-1}(t)})\}|=:
I_1+I_2.\end{split}\end{equation} By (\ref{L4}) and $
\xi_y^{-1}(t)\le t$ we obtain

\beq\label{L7} I_1\le \|\nn f\|_\infty \e^{K' t}R
\rr(x,y).\end{equation} Moreover, since $f\in C_0^\infty(M)$ with
$Nf|_{\pp M}=0,$ it follows from the It\^o formula and (\ref{L5})
that

$$I_2\le \bigg|\EE\int_{\xi_x^{-1}(t)\land
\xi_y^{-1}(t)}^{\xi_x^{-1}(t)\lor \xi_y^{-1}(t)} L' f(X_s')\d
s\bigg|\le \|L'f\|_\infty \EE|\xi_x^{-1}(t)-\xi_y^{-1}(t)|\le
c_1t\e^{K't} \rr(x,y)$$ holds for some constant $c_1>0.$  Combining
this with (\ref{L6}) and (\ref{L7}) we conclude that

$$\|\nn P_tf\|_\infty\le c_2 (1+t) \e^{K't},\ \ \ t\ge 0$$ for some
constant $c_2>0$.\end{proof}

\ \newline\emph{Proof of Prposition \ref{A2}.}  Let $f\in C_b^1(M)$. By Lemma \ref{L2.1}
 we only have to prove the  boundedness
of $|\nn P_\cdot f|$ on $[0,T]\times M$ .

 (a) Let $f\in C_0^\infty(M).$ In this case there exist a sequence of functions
  $\{f_n\}_{n\ge 1}\subset C_0^\infty(M)$ such that $Nf_n|_{\pp M}=0, f_n\to f$ uniformly as $n\to\infty$,
  and $\|\nn f_n\|_\infty\le 1+\|\nn f\|_\infty$ holds
  for any $n\ge 1,$ see e.g. \cite{W94}. By Lemmas \ref{L2.1} and \ref{L2.3}, (\ref{Hsu})
  holds for $f_n$ in place of $f$ so that Proposition \ref{AA1} implies

  $$\ff{|P_t f_n(x)-P_t f_n(y)|}{\rr(x,y)}\le C,\ \ \ t\le T, n\ge 1, x\ne y$$ for some constant $C>0.$
  Letting first $n\to 0$ then $y\to x$, we conclude that $|\nn P_\cdot f|$ is bounded
  on $[0,T]\times M.$

  (b) Let $f\in C_b^\infty(M)$.  Let $\{g_n\}_{n\ge 1}\subset C_0^\infty(M)$ be such that
  $0\le g_n\le 1, |\nn g_n|\le 2$ and $g_n\uparrow 1$ as $n\uparrow\infty$. By (a) and Lemma \ref{L2.1},
  we may apply (\ref{Hsu}) to $g_nf$ in place of $f$
  such that Proposition  \ref{AA1} implies

   $$\ff{|P_t (g_nf)(x)-P_t(g_n f)(y)|}{\rr(x,y)}\le C,\ \ \ t\le T, n\ge 1, x\ne y$$
   holds for some constant $C>0.$ By the same reason as in (a) we conclude that $|\nn P_\cdot f|$ is bounded
  on $[0,T]\times M.$

(c) Finally, for $f\in C_b^1(M)$ there exist $\{f_n\}_{n\ge
1}\subset C_b^\infty(M)$ such that $f_n\to f$ uniformly as
$n\to\infty$ and $\|\nn f_n\|_\infty\le \|\nn f\|_\infty+1$ for any
$n\ge 1.$ Therefore, the proof is complete by the same reason as in
(a) and (b). \qed

\beg{thebibliography}{99}

\bibitem{B} D. Bakry, \emph{On Sobolev and logarithmic Sobolev inequalities for Markov semigroups, }
New Trends in Stochastic Analysis, 43Ð75, World Scientific, 1997.

\bibitem{BE} D. Bakry and M. Emery, \emph{Hypercontractivit\'e de
semi-groupes de diffusion}, C. R. Acad. Sci. Paris. S\'er. I Math.
299(1984), 775--778.

\bibitem{BL} D. Bakry and M. Ledoux, \emph{L\'evy-Gromov's isoperimetric inequality for an
infinite-dimensional diffusion generator,} Invention Math. 123(1996), 259--281.

\bibitem{BGL} S. G. Bobkov, I.  Gentil and M.  Ledoux,
\emph{Hypercontractivity of Hamilton-Jacobi equations,} J. Math.
Pures Appl. 80(2001), 669--696.

\bibitem{Hsu} E. P. Hsu, \emph{Multiplicative functional for the heat
equation on manifolds with boundary,} Michigan Math. J. 50(2002),
351--367.

 \bibitem{K2} A. Kasue, \emph{Applications of Laplacian and Hessian comparison theorems,} Geometry of
 geodesics and related topics (Tokyo, 1982), 333--386, Adv. Stud. Pure Math., 3, North-Holland, Amsterdam, 1984.

 \bibitem{K2} A. Kasue, \emph{On Laplacian and Hessian comparison theorems,}
Proc. Japan Acad. 58(1982), 25--28.

\bibitem{Ledoux} M. Ledoux, \emph{The geometry of Markov diffusion generators,} Ann. Facu. Sci. Toulouse
9(2000), 305--366.

 \bibitem{OV} F.  Otto and C. Villani, \emph{Generalization of an inequality
by Talagrand and links with the logarithmic Sobolev inequality,} J.
Funct. Anal.  173(2000), 361--400.
\bibitem{Q} Z. Qian, \emph{A gradient estimate on a manifold with convex boundary,}
 Proc. Roy. Soc. Edinburgh Sect. A 127(1997), 171--179.

\bibitem{S} M.-K. von Renesse and K.-T.  Sturm, \emph{Transport inequalities, gradient estimates, entropy,
and Ricci curvature,}  Comm. Pure Appl. Math. 58 (2005),  923--940.

\bibitem{W94} F.-Y. Wang, \emph{Application of coupling methods to the Neumann eigenvalue
  problem}, Probab. Theory Related Fields 98 (1994),  299--306.

\bibitem{W97} F.-Y. Wang,  \emph{Logarithmic {S}obolev inequalities on
noncompact
  {R}iemannian manifolds}, Probab. Theory Related Fields 109(1997), 417--424.

\bibitem{W04} F.-Y. Wang, \emph{Equivalence of dimension-free Harnack inequality and curvature condition,}
 Integral Equations Operator Theory 48 (2004),  547--552.

 \bibitem{W05}  F.-Y. Wang, \emph{Gradient estimates and the first
Neumann eigenvalue on manifolds with boundary,} Stoch. Proc. Appl.
115(2005), 1475--1486.

 \bibitem{W07} F.-Y. Wang, \emph{Estimates of the first Neumann eigenvalue and the log-Sobolev constant
  on  nonconvex manifolds,} Math. Nachr. 280(2007), 1431--1439.

\bibitem{W09} F.-Y. Wang, \emph{Second fundamental form and gradient of Neumann semigroup,}  J. Funct. Anal.
256(2009), 3461--3469.


\bibitem{W09b} F.-Y. Wang, \emph{Semigroup properties for the second fundamental form,} arXiv:0908.2890.

\end{thebibliography}

\end{document}